\documentclass[12pt]{amsart}
\usepackage[dvips]{graphicx}
\usepackage{amssymb}
\def\C{\mathbb C}

\def\N{\mathbb N}

\def\R{\mathbb R}
\usepackage{psfig}
\usepackage{latexsym}
\usepackage{amsfonts}
\usepackage{amssymb}

\def\dist{\mathrm{dist\,}}
\def\diam{{\mathrm{diam}\,}}

\def\bC{{\overline{\C}}}
\def\Re{\operatorname{Re}}
\def\Im{\operatorname{Im}}
\def\:{\colon}

\newtheorem{lem}{Lemma}
\theoremstyle{remark}
\newtheorem*{ack}{Acknowledgment}
\title{Meromorphic functions with two completely invariant
domains}
\dedicatory{Dedicated to the memory of Professor I.~N.~Baker}
\begin{document}
\begin{abstract}
We show that if a meromorphic function has two completely 
invariant Fatou components  and only finitely many 
critical and asymptotic
values, then its Julia set is a Jordan curve.
However, even if both domains are attracting basins,
the Julia set need not be a quasicircle.
We also show that
all critical and asymptotic values are contained in
the two completely invariant components.
This need not be the case for functions with infinitely 
many critical and asymptotic values.
\end{abstract}
\author{Walter Bergweiler}\thanks{Supported  by the German-Israeli Foundation
for Scientific Research and Development (G.I.F.),
grant no.\ G -643-117.6/1999.}
\address{Mathematisches Seminar,
Christian--Albrechts--Universit\"at zu Kiel,
Ludewig--Meyn--Str.~4,
D--24098 Kiel,
Germany}
\email{bergweiler@math.uni-kiel.de}
\author{Alexandre Eremenko}\thanks{
Supported by NSF grants DMS-0100512
and DMS-0244421, and by the Humboldt Foundation.}
\address{Department of Mathematics,
Purdue University, West Lafayette, IN 47907, USA}
\email{eremenko@math.purdue.edu}
\maketitle

\section{Introduction and main result}
Let $f$ be a meromorphic function in the complex plane $\C$.
We always assume that $f$ is not fractional linear or constant.
For the definitions and main facts of the theory of
iteration of meromorphic functions we refer to 
a series of papers by Baker, Kotus and L\"u
\cite{Bak91,Bak90b,BKL,Bak92},
who started the subject,
and to the survey article~\cite{Bergweiler}.
For the dynamics of rational functions we refer to 
the books
\cite{Bea91,CG,Milnor, Stein93}.

A {\em completely invariant domain}
is a component $D$ of the set of normality such that
$f^{-1}(D)=D$.
There is an unproved conjecture (see \cite[p.~608]{BKL},
\cite[Question~6]{Bergweiler}) that a meromorphic function can
have at most two completely invariant domains.
For rational functions this fact easily follows from
Fatou's investigations \cite{Fatou},
and it was first explicitly
stated by Brolin \cite[\S 8]{Brolin}.
Moreover, if a rational function has two completely invariant
domains, then their common boundary is a Jordan curve on the Riemann
sphere,
and each of the domains coincides with the basin of attraction
of an attracting or superattracting fixed point, or of an
attracting petal of a neutral fixed point with multiplier $1$;
see \cite[p. 300-303]{Fatou} and \cite{Brolin}.
All critical values of
$f$ are contained in the completely invariant domains.

In this paper we extend these results to a class
of transcendental meromorphic functions
in $\C$. 
This class $S$ consists
of meromorphic functions with finitely many critical and
asymptotic values. Let $A=A(f)$ be the set
of critical and asymptotic values. 
We also call the elements of $A$ {\em singular values} of $f$.
For $f\in S$ the map
$$f:\C\backslash f^{-1}(A)\to\bC\backslash A$$
is a covering. 
By $J=J(f)\subset \C$ we denote the Julia set of $f$.

Baker, Kotus and L\"u \cite[Theorem~4.5]{BKL} proved 
that functions of the class $S$ have
at most two completely invariant domains.
Cao and Wang~\cite[Theorem~1]{CW} have shown that if a function in the 
class $S$ has two completely invariant domains, then its
Fatou set is the union of these domains.
We complement these results with the following

\medskip

\noindent
{\bf Theorem.} {\em
Let $f$ be a function of the class $S$,
having two completely invariant domains $D_j,\; j=1,2.$
Then
\begin{itemize}
\item[$(i)$] 
each $D_j$ is the basin
of attraction of an attracting or
superattracting fixed point, or of a petal of a neutral
fixed point with multiplier $1$,
\item[$(ii)$] 
$A(f)\subset D_1\cup D_2$,
\item[$(iii)$] 
each $D_j$ contains at most one asymptotic value,
and if
$a$ is an asymptotic value and 
$0<\epsilon< \dist(a,A\backslash\{a\})$,
then the set
$\{ z:|f(z)-a|<\epsilon\}$ has only one unbounded component,
\item[$(iv)$] 
$J\cup\{\infty\}$ is a Jordan curve in $\bC$.
\end{itemize}
}

\medskip

A simple example of a meromorphic function
of class $S$ with two
completely invariant domains is $f(z)=\tan z$, for which
the upper and lower half-planes are
completely invariant, and
each of these half-planes is attracted to one of 
the two petals of the fixed point $z=0$.
 
More examples will be given later in \S \ref{examples}.

In the case that $f$ is rational and both $D_1$ and $D_2$
are attracting or superattracting basins, 
Sullivan \cite[Theorem~7]{Sul}
and 
Yakobson \cite{Yakob}
proved that $J$ is a quasicircle.
Steinmetz \cite{Stein}
extended this result to the case that both
completely invariant domains are basins of
two petals attached to the same neutral fixed point.
We will construct
examples of transcendental functions  in $S$
for which  $D_1$ and $D_2$ 
are attracting basins, or basins of petals attached to
the same neutral fixed point, but where $J$ is not 
a quasicircle; see Examples 1 and 2 in \S \ref{examples}.

On the other hand,
Keen and Kotus \cite[Corollary~8.2]{KK}
have shown that for the family $f_\lambda(z)=\lambda
\tan z$ there exists a domain $\Omega$ containing 
$(1,\infty)$ such that 
$f_\lambda$ has two completely invariant 
attracting basins and
$J(f_\lambda)$ is a quasicircle
for $\lambda\in\Omega$.
Meromorphic functions for which the Julia set is contained
in a quasicircle were also considered by Baker, Kotus and
L\"u~\cite[\S 5]{Bak91}.

Baker \cite{Bak70a} proved that an entire function $f$
can have at most one completely invariant component of the set
of normality,
and that
such a domain contains all critical values.
Eremenko and Lyubich
\cite[\S 6]{EL}
proved that a completely invariant
domain of an entire function
also contains all asymptotic values of a certain type,
namely
those associated with direct singularities of $f^{-1}$.
On the other hand, Bergweiler \cite{Bergweiler2}
constructed
an entire function with a completely invariant
domain $D$, and
such that some asymptotic value
belongs to the Julia set $J=\partial D$.
Example~ 3 in \S 3 shows that
meromorphic functions with two completely
invariant
components of the set of normality
can have asymptotic values on their
Julia sets.
So (ii) does not hold for general
meromorphic functions,
without the assumption that $f\in S$.

\section{Proof of the Theorem}
We shall need the following result of Baker, Kotus 
and L\"u~\cite[Lemmas 4.2 and 4.3]{BKL}
which does not require that $f\in S$.
Here and in the following
all topological
notions are related to $\C$ unless $\overline{\C}$
is explicitly mentioned.
\begin{lem}
\label{lemma1}
Let $f$ be meromorphic with two completely invariant components
$D_1$ and $D_2$ of the set of normality.
Then $D_1$ and $D_2$ are simply-connected and
$J=\partial D_1=\partial D_2$. In particular, $J$
is a connected subset of $\C$.

\end{lem}

{\em Proof of the Theorem}.
As the result is known for rational $f$,
we assume that our $f$ is transcendental.

Statement
(i) follows from the classification
of dynamics on the Fatou
set for meromorphic functions of the class $S$,
given in
\cite[Theorems 2.2 and 2.3]{Bak92},
\cite[Theorem 6]{Bergweiler},
and \cite[p.~3252]{RS}.

To prove (ii), we consider for $j=1,2$ the finite sets
$A_j=A\cap D_j$. Let $\Gamma_j$
be a Jordan curve in $D_j$ which separates
$A_j$ from $\partial D_j$. Let $G_j$ be the Jordan
regions bounded by the $\Gamma_j$. Let
$G=\overline{\C}\backslash
(\overline{G_1}\cup\overline{G_2})$
be the doubly connected
region in $\overline{\C}$ bounded by $\Gamma_1$
and $\Gamma_2$.
Notice that $G$ contains the Julia set $J$.

We define $\gamma_j=f^{-1}(\Gamma_j)$. 
Then
\begin{equation}\label{1}
f:\gamma_j\to \Gamma_j,\; j=1,2,
\end{equation}
are covering maps.

We claim that each $\gamma_j\subset\C$
is a single simple curve
tending to infinity in both directions,
which means that $\gamma_j\cup\{\infty\}$
is a Jordan curve in $\bC$.

To prove the claim, we fix $j$ and consider
the full preimages
$H_j=f^{-1}(G_j)$ and
$F_j=f^{-1}(D_j\backslash \overline{G_j})$.
Then $D_j=F_j\cup H_j\cup\gamma_j$. The boundary of each component of
$F_j$ contains a component of $\gamma_j$, and this gives a bijective
correspondence between components of $F_j$
and components of $\gamma_j$.

We notice that
$H_j$ is connected. Indeed,
by complete invariance of $D_j$, we have
$H_j\subset D_j$, so
every two points $z_1$ and $z_2$
in $H_j$ can be connected by a curve $\beta$ in $D_j$,
so that $\beta$
does not pass through the critical points of $f$.
The image  $f(\beta)$ of this curve begins and ends
in $G_j$, and
does not pass through the critical values of $f$.
By a small perturbation of
$\beta$ we achieve that $f(\beta)$ does not
pass through asymptotic values.
Using the fact that
$$f:D_j\backslash f^{-1}(A_j)\to D_j\backslash A_j,$$
is a covering
and that $A_j\subset G_j$, we can deform $\beta$ 
into a curve
in $H_j$ which still connects $z_1$ and $z_2$.
This proves that the $H_j$ are connected.

It follows that $H_j$ is unbounded,
as it contains infinitely many
preimages of a generic point in $G_j$.

It is easy to see that the boundary of each
component $F_j^\prime$ of $F_j$ intersects the
Julia set.

For each component $F_j^\prime$ of
$F_j$, the intersection $\partial F_j^\prime\cap \overline{H_j}$
is a component $\gamma_j^\prime$ of $\gamma_j$.
This component $\gamma_j^\prime$ divides the plane into two
parts, one containing $H_j$ and the other containing
$F_j^\prime$.
We conclude that every component of $\gamma_j$
is unbounded,
because $H_j$ is unbounded,
and $\partial F_j$ intersects the Julia
set which is unbounded and connected by Lemma 1.
(A similar argument for unboundedness of each
component of $\gamma_j$ was given in \cite{BKL}).  

For every component $\gamma_j^\prime$ of $\gamma_j$,
the component of $\C\backslash\gamma_j^\prime$
that contains $F_j^\prime$ intersects the Julia set.
Since the Julia set is connected by Lemma~1,
we conclude that $F_j$ and
$\gamma_j$ are connected.
So the map (\ref{1}) is
a universal covering by a connected set
$\gamma_j$, for each $j=1,2$. 
Thus $\gamma_1\cup\{\infty\}$ and
$\gamma_2\cup\{\infty\}$
are Jordan curves in $\bC$ whose intersection
consists of the single point $\infty$. 
This proves our claim.

As a corollary we obtain that the point $\infty$
is accessible from each $D_j$, and so all
poles of all iterates $f^{n}$ are accessible from
each $D_j$. (This fact was established in~\cite{BKL}.)

Next we note that
 the set
$\gamma_1\cup\gamma_2\cup\{\infty\}$ 
separates the sphere into three simply connected
regions. We denote by $W$ that region whose
boundary in $\bC$
is $\gamma_1\cup\gamma_2\cup\{\infty\}$.
Then
\begin{equation}
\label{3}
f^{-1}(G)=W,
\end{equation}
in particular, $W$ contains the Julia set $J$.

To complete the proof of (ii) we choose an arbitrary point
$w\in J$ and show that $w$
is neither a critical value nor an asymptotic value.

Fix an arbitrary point $w_1\in\Gamma_1$.
The preimage 
$f^{-1}(w_1)$ consists of infinitely
many points $a_k\in\gamma_1$,
which we enumerate by all integers
in a natural order on $\gamma_1$. Let $\phi_k$ be
the branches of $f^{-1}$ such that $\phi_k(w_1)=a_k$.
We find a simple
curve $\Delta$ from $w_1$ to some point
$w_2\in \Gamma_2$ such that
$\Delta\backslash\{ w_1,w_2\}$
is contained in
$G\backslash\{ w\}$,
and such that
all branches $\phi_k$ have analytic
continuation along $\Delta$
to the point $w_2$. We define
$$G'=G\backslash\Delta\subset\bC.$$
The full preimage $f^{-1}(\Delta)$
consists of infinitely
many disjoint simple curves $\delta_k$
starting at the points
$a_k$ and ending at some points $b_k\in\gamma_2$.
The open curves $\delta_k\backslash\{ a_k,b_k\}$ are
disjoint from $\gamma_1\cup\gamma_2$.

For every integer $k$, let $Q_k$ be the Jordan region
bounded by $\delta_k,\,\delta_{k+1}$, the arc $(a_k,a_{k+1})$
of $\gamma_1$ and the arc
$(b_k,b_{k+1})$ of $\gamma_2$.
Then $f$ maps $Q_k$ into $G'$,
and $f(\partial Q_k)\subset
\partial G'$. So 
\begin{equation}\label{2}
f:Q_k\to G'
\end{equation}
is a ramified covering, continuous up to the boundary.
Furthermore, the boundary map is a local homeomorphism.
As each point of $\Gamma_1\backslash\{ w_1\}$ has only
one preimage on $\partial Q_k$, we conclude that
(\ref{2}) is a homeomorphism.
Now it follows that the restriction
$f:(b_k,b_{k+1})\to\Gamma_2\backslash\{ w_2\}$
is a homeomorphism and
$$W=\bigcup_{k=-\infty}^\infty Q_k\cup \delta_k\backslash\{a_k,b_k\}.$$
It follows that there are no critical points over $w$,
so $w$ is not a critical value.

If $w$ were an asymptotic value,
there would be a curve
$\alpha$ in $W$
which tends to infinity, and such that $f(z)\to w$ as
$z\to\infty,\; z\in\alpha$.
But this curve $\alpha$ would
intersect infinitely many of
the curves $\delta_k$, so its image $f(\alpha)$
would intersect
$\Delta$ infinitely many times,
which contradicts the assumption
that $f(\alpha)$ tends to $w$.

This completes the proof of (ii).
The proof actually shows that
$f:W\to G$ is a universal covering,
a fact which we will use
later.

To prove (iii), let us assume that 
$D_1$ contains two asymptotic values, or that
$\{z\in D_1:|f(z)-a|<\epsilon\}$ has 
two unbounded components for some asymptotic 
value $a\in D_1$.
Then there
exists a curve $\alpha\subset D_1$, tending to infinity
in both directions, such that $f(z)$ has
limits as 
$z\to\infty$,
$z\in\alpha$,
 in both directions,
where these limits are the two asymptotic values 
in the first case, and where both limits are equal to
$a$ in the second case, but the two tails of the curve
$\alpha$ are in different components of 
$\{z\in D_1:|f(z)-a|<\epsilon\}$.
Now one of the regions, say $R$,
into which $\alpha$ partitions
the plane does not intersect the Julia set $J$
(because $J$ is connected by Lemma~\ref{lemma1}),
and thus $R\subset D_1$.
We want to conclude that
$f$ has a limit as $z\to\infty$ in $R$.

To do this, we choose an arbitrary point $b\in D_2$
and consider the function $g(z)=(f-b)^{-1}$ which
is holomorphic and bounded in $D_1$.
This function has limits when $z\to\infty$, $z\in\alpha$,
so by a theorem of Lindel\"of \cite{Nev}, these limits coincide
and $g$ has a limit as $z\to\infty$ in $R$.
This proves  (iii).

To prove (iv), we distinguish several cases, according
to the dynamics of $f$ in each $D_j$.

\medskip

\noindent
{\bf 1}.
Suppose first that both $D_1$ and $D_2$
are basins of attraction
of attracting or superattracting points.
Then we choose the curves $\Gamma_j$ as above,
but with the 
additional property that $f(\Gamma_j)\subset G_j$,
that is $f(\Gamma_j)\cap\overline{G}=\emptyset.$
To achieve this, we denote by $z_j$
the attracting or superattracting fixed point in $D_j$,
choose $G_j$ to be the
open hyperbolic disc centered
at $z_j$, of large enough hyperbolic radius, so that 
$A_j\subset G_j$, and put $\Gamma_j=\partial G_j$.
Then the $G_j$ are $f$-invariant,
and moreover $f(\overline{G_j})
\subset G_j$
for $j=1,2$, because $f$ is strictly contracting
the hyperbolic metric in $D_j$. 
It follows that the closure of $W=f^{-1}(G)$
is contained in $G$.
Let $h$ be the 
hyperbolic metric
in $G$, and $|f'(z)|_h$ the infinitesimal length
distortion by $f$ at the
point $z\in W$ with
respect to $h$.
By the Theorem of Pick
\cite[Theorem 2.11]{Milnor}
there exists $K>1$ such that
\begin{equation}
\label{4}
|f'(z)|_h\geq K,\quad z\in W.
\end{equation}
Now we consider successive preimages $W_n=f^{-n}(W)$.
Note that $\infty\notin A(f)$ by (ii) which implies that the 
components of $f^{-1}(\gamma_j)$ are bounded for $j=1,2$.
We deduce that every component of $W_n$ is a Jordan domain whose
boundary consists of two
cross-cuts, one of $D_1$ and the other of $D_2$.
These crosscuts meet at two poles of $f^n$.
It follows from (\ref{4}) that the diameter
(with respect to the metric
$h$) of every component
of $W_n$ is at most $CK^{-n}$, where $C>0$ is a constant.
Now we notice that 
$$J=\bigcap_{n=1}^\infty\overline{W_n}$$
and prove that every point $z\in J$ is accessible
both from $D_1$ and $D_2$.

The accessibility of poles of the iterates $f^{n}$
was already noticed before.
Now we assume that $z$ is not a pole of any iterate.
Let $V_n$ be the component of $W_n$ that contains $z$.
Then
$V_1\supset V_2\supset\ldots$.
The intersection $V_k\cap D_j$ is connected (its
relative boundary with respect to $D_j$
is a cross-cut in $D_j$),
so one can choose a sequence $z_{k,j}\in V_k\cap D_j$
and connect
$z_{k,j}$ with $z_{k+1,j}$ by a curve $\ell_{k,j}$
in $V_k\cap D_j$. The union of these curves gives
a curve in $D_j$ which tends to $z$.

The proof of (iii) in the attracting case
is completed by an application of Schoenflies' theorem
that if each point of a
common boundary of two domains
on the sphere is accessible
from both domains then this common boundary is
a Jordan curve \cite{New}.

\medskip

\noindent
{\bf 2}. To prove (iv) in the remaining cases, 
suppose, for example, that $D_1$
is the domain of attraction of a petal associated with
a neutral fixed point $a$. We need several lemmas.

\begin{lem}
There exists a Jordan domain $G_1$
with the properties $\overline{G_1}\subset D_1\cup\{ a\}$,
$f(\overline{G_1})\subset G_1\cup\{ a\}$, $A_1\subset G_1$,
and $G_1$ is absorbing,
that is for every compact $K\subset D_1$ there exists a
positive integer $n$ such that $f^{n}(K)\subset G_1$.
\end{lem}

{\em Proof}.
It is well known (see, e.~g.~\cite[\S 10]{Milnor})
that there exists a domain $G_1$ having all properties
mentioned except possibly $A_1\subset G_1$. 
Such a domain is called an {\em attracting petal}.

Let $P$ be an attracting petal.
Choose a point $z_0\in D_1$ and let
$r>0$ be so large that the open hyperbolic
disc $B(z_0,r)$
of radius $r$
centered at $z_0$
contains $A_1$, and $z_1=f(z_0)\in B(z_0,r)$.
Then put 
$$G_1^\prime=P\cup\left(\bigcup_{k=0}^\infty B(z_k,r)\right),
\quad z_k=f^{k}(z_0).$$
Then $G_1^\prime$ is absorbing because the petal $P$ is absorbing.
Notice that for every neighborhood
$V$ of $a$, all but finitely
many discs $B(z_k,r)$ are
contained in $V$.
This easily follows
from the comparison of the Euclidean
and hyperbolic metrics near
$a$, or, alternatively, from the local description
of dynamics near a neutral fixed point with multiplier
$1$.
It is easy to see that
$f(\overline{G_1^\prime})\subset G_1^\prime
\cup\{ a\}$. 
Now we fill the holes in $G_1^\prime$: let $X$ be the unbounded
component of $\overline{G_1^\prime}$ and
$G_1=\C\backslash\overline{X}$.
It is easy to see that $G_1$ is a Jordan domain
(its boundary is
a union of arcs of hyperbolic circles which is locally finite,
except at the point $a$, plus some boundary arcs of the
petal).~\hfill$\Box$

\medskip

Now we fix the following notations till the
end of the proof of the Theorem. If $D_j$ is a
basin of an attracting or superattracting fixed point,
let $G_j$ be the Jordan region constructed
in the first part of the proof of (iv).
If $D_j$ is a basin of a petal, let $G_j$
be the region from Lemma 2.
We define $\Gamma_j=\partial G_j$. This is a
Jordan curve in $D_j$ or in 
$D_j\cup\{ a\}$ which encloses all
singular values in $D_j$.

Next we define $G=\bC\backslash (\overline{G_1}\cup
\overline{G_2})$. This region is simply connected
in the case that both $D_1$ and $D_2$ are basins
of two petals associated with the {\em same} fixed
point, and doubly connected in all other cases.
If $G$ is doubly connected, we make a simple cut
$\delta$ disjoint from the set $A$ of
singular points, as in the proof of (ii),
to obtain a simply
connected region $G'=G\backslash\delta$.
If $G$ is simply connected we set $G'=G$.
All branches of $f^{-n}$ are holomorphic in $G'$.
Let $\gamma_j=f^{-1}(\Gamma_j)$.

\begin{lem}
There exists a repelling fixed point
$b\in J$ which is accessible from both $D_1$ and $D_2$
by simple 
curves $\beta_j$ 
which begin at
some points of $\Gamma_j$ and do not intersect
$G_j$,
and which satisfy
$f(\beta_j)\cap \overline{G'}=\beta_j$,
 for $j=1,2$.
\end{lem}

{\em Proof}. We use the notation introduced before the
statement of the Lemma. 
Fix one of the components $Q$, of $f^{-1}(G')$,
such that $\overline{Q}\subset G'$. 
Let $\phi$ be the branch of $f^{-1}$ which maps
$G'$ onto $Q$. Then $\phi$ has an attracting fixed
point $b\in Q$. Let $z_0\in \Gamma_1$ and
$z_1=\phi(z_0)\in \partial Q$.
We connect $z_0$ and $z_1$ by a simple curve
$\beta$ in $(G'\backslash Q)\cap D_1$. Such a
curve exists because $z_1\in\gamma_1$,
and the component of $\C\backslash\gamma_1$ that contains
$G_1$ is completely contained in $D_1$.

Now $$\beta_1=\bigcup_{k=1}^\infty \phi^k(\beta)$$
is a curve in $D_1$ tending to $b$ which 
satisfies $f(\beta_1)\cap \overline{G'}=\beta_1$.
Similarly a curve $\beta_2$ in $D_2$ is constructed.~\hfill$\Box$

\medskip

\def\Gpp{{G^{\prime\prime}}}
Now, if $G$ is doubly connected,
we set $\Gpp=G\backslash(\beta_1\cup\beta_2\cup
\{ b\}).$ If $G$ is simply connected then $\Gpp=G$.
Then $\Gpp$ is a simply connected region which contains
no singular values of $f$. Let $\{\phi_k\}_{k\in \N}$
be the set of all branches of $f^{-1}$ in $\Gpp$. These
branches map $\Gpp$ onto Jordan regions $T_k\subset
\Gpp$. These regions $T_k$ are of two types:
the regions of the first type are contained in $\Gpp$ with
their closures, while the regions of the second type
have common boundary points with $\Gpp$.

We claim that there are only finitely many regions of
the second type. To study these regions $T_k$, we first observe
that the full preimage of $\Gamma_j$ is a curve $\gamma_j$
which can have at most one point in common with $\Gamma_j$,
namely the neutral fixed point on $\Gamma_j$.
Thus the region $W=f^{-1}(G)$ bounded by $\gamma_1$ and $\gamma_2$ is
a simply connected region contained in $G$, and the boundary
$\partial W$ has at most two common points with $\partial G$,
namely the neutral fixed points. The full preimage of the
cross-cut 
$$\alpha=\beta_1\cup\beta_2\cup\{ b\}$$
constructed in Lemma 2 consists of
countably many disjoint curves $\alpha_k\subset\overline{W}$.
Each $\alpha_k$ connects a point on $\gamma_1$ to a point on
$\gamma_2$. One of the $\alpha_k$, say $\alpha_1$, is contained
in $\alpha$ while all others are disjoint from $\alpha$.
Thus our regions $T_k$ are curvilinear rectangles,
similar to the $Q_k$ used in the proof of (ii). 
In particular, they cluster only at $\infty$ so that only
finitely many of them are of the first type.

It is easy to see that every region of 
the second type has on
its boundary exactly one of the following points:
a neutral fixed point or the repelling fixed point $b$.
Indeed, let $T$ be a region of the second type,
and $\phi:\Gpp\to T$ the corresponding branch of the
inverse. Then the iterates $\phi^n(z)$ converge to
a unique point $c\in \overline{T}$ by the Denjoy--Wolff Theorem.
(This theorem is usually stated for the unit disk,
but it follows for Jordan domains like $T$ by the 
Riemann mapping theorem, using that the Riemann map
extends homeomorphically to the boundary.)
On the other hand, it follows from the local dynamics
near the repelling fixed point $b$ and a  neutral fixed point $a$
that there exists $\epsilon>0$ such
that 
$\phi^n(z)\to b$ if $|z-b|<\epsilon$ and
$\phi^n(z)\to a$ if $|z-a|<\epsilon$, $z\in T$.

If $\phi_j$ and $\phi_k$ are two different branches of $f^{-1}$ in $\Gpp$,
whose images are of the second type, then the 
images $(\phi_k\circ\phi_j)(\Gpp)$ are compactly contained in
$\Gpp$. There exists a compact subset set $K\subset\Gpp$ which contains
all regions $T$ of the first type as well as all images
$(\phi_k\circ\phi_j)(\Gpp)$ where $j\neq k$.

Now consider the hyperbolic metric in $\Gpp$ and let $|\phi'|_h$ stand
for the infinitesimal length distortion of a branch $\phi$ with respect to
this hyperbolic metric. Then we have for all $z\in\Gpp$ and some
$\lambda\in(0,1)$:
\begin{equation}\label{dist1}
|\phi_k^\prime(z)|_h<\lambda\quad\mbox{for all $k$ of the first type}
\end{equation}
and
\begin{equation}\label{dist2}
|(\phi_j\circ\phi_k)'(z)|_h<\lambda\quad\mbox{for all $j\neq k$}.
\end{equation}

Let $W_n=f^{-n}(W)$. Then the Julia set 
can be represented as the intersection of a decreasing sequence of
closed sets $J=\bigcap_{n=1}^\infty \overline{W_n}$. The points of the
Julia set are divided into the following categories:

a) poles of $f$ and their preimages,

b) neutral fixed points and their preimages, 

c) the repelling point $b$ and its preimages

d) those points of $J$ which are interior to all $f^{-n}(\Gpp)$.

We have already seen that all points of the categories
a)-c) are accessible from each of the domains $D_1$ and $D_2$.

The proof that the points of the type d) are accessible is similar to
the argument in the case that both $D_1$ and $D_2$ are
attracting basins: we will show that  each such point $z$ can be surrounded
by a nested sequence of Jordan curves whose diameter tends to zero.

Indeed, each point $z$ of the class d) can be obtained as a limit
$$z=\lim_{n\to\infty}(\phi_{k_1}\circ\phi_{k_2}\circ\ldots
\circ\phi_{k_n})(w),$$ 
where $w\in\Gpp$. For a point $z$ of the
category d), the sequence $k_1,k_2,\ldots$ is uniquely defined.
We will call this sequence the {\em itinerary} of $z$.
Let us consider the domains
$$T_n(z)=(\phi_{k_1}\circ\phi_{k_2}\circ\ldots\circ\phi_{k_n})(\Gpp),$$
in other words, $T_n(z)$ is that component of $f^{-n}(\Gpp)$ which
contains $z$. The boundary of $T_n(z)$ is a Jordan curve which
intersects the Julia set at a finite set of points of categories
a)-c). The complementary arcs of these points are cross-cuts
of $D_1$ and $D_2$. Thus, to show that $z$ is accessible
from $D_1$ and $D_2$, it is enough to show that 
the diameter
of $T_n(z)$ tends to zero as $n\to\infty$. Let $z\in J$
be a point of
category d), and $k_1,k_2,
\ldots$ its itinerary.
Then the sequence $k_1,k_2,\ldots$ cannot have
an infinite tail consisting of the branch numbers of 
the second type.
Indeed, the iterates of any branch of 
the second type converge
to a boundary point $x$ of $\Gpp$ (a neutral fixed point or the
point $b$). In this case, $z$ will be a preimage of $x$.

Since the itinerary does not stabilize on a branch number
of the second type, we can use (\ref{dist1}) and (\ref{dist2}) to conclude
that $\diam T_n(z)\to 0$.

This completes the proof.~\hfill$\Box$

\section{Examples}
\label{examples}

\noindent
{\bf Example 1.}
Let
$$g(z)=\frac{1}{1+a\cos\sqrt{z}}$$
where $0<a<\frac15$.
Then there exists $b<0$ such that
$$f(z)=\frac{g(z+b)-g(b)}{g'(b)}$$
has a parabolic fixed point at zero, with two completely invariant
parabolic basins attached to it.
Moreover, $f\in S$ and the Julia set of $f$ is a Jordan curve,
but not a quasicircle.

\medskip
{\em Verification.}
Note that
$g$ has no poles on the real axis.
We have
$$g'(z)=\frac{a\sin\sqrt{z}}{\sqrt{z}(1+a\cos\sqrt{z})^2}$$
and
$$g''(z)=-\frac{a^2(\cos\sqrt{z})^2-a\cos\sqrt{z}-2a^2}{z(1+a\cos\sqrt
{z})^3}
-\frac{a\sin\sqrt{z}}{z\sqrt{z}(1+a\cos\sqrt{z})^2}.$$
It follows that
$$\lim_{x\to -\infty} g''(x) x\cos\sqrt{x}=-\frac{1}{4a}<0$$
so that $g''(x)>0$ if $x$ is negative and of sufficiently
large modulus.
On the other  hand,
$$g''(0)=\frac{a(5a-1)}{2(a+1)^3 }<0.$$
Thus there exists $b\in(-\infty,0)$ with
$g''(b)$=0 and
$g''(x)>0$ for $x<b$.

The critical points of $g$ are given by
$(k\pi)^2$ where $k\in\N$,
and $g$ has a maximum there for odd $k$ and a
minimum for even $k$.
It follows that $g'(x)>0$ for $x<\pi^2$
and thus in particular for $x\leq b$.
Thus
$f$ has the
critical points $(k\pi)^2-b$, with corresponding
critical values
$$d_\pm=\frac{(1\pm a)^{-1}-g(b)}{g'(b)}.$$
Moreover, $f$ has the asymptotic value
$c=-g(b)/g'(b)$,
which is also a Picard exceptional value of $f$,
and no other asymptotic values.
Thus $A(f)=\{c,d_+,d_-\}$.

Next we note that $f(0)=0$, $f'(0)=1$ and $f''(0)=0$.
Since $f'(x)=g'(x+b)/g'(b)$ we have $0<f'(x)<1$ for $x<0$.
It follows from the mean value theorem
that if $x<0$, then
$x<f(x)<0$.
Thus $(-\infty,0)$ lies in a parabolic basin $U$ attached to
the parabolic point $b$.
In particular, $U$ contains the value
$c=-g(b)/g'(b)$ which is a Picard exceptional value of $f$.
We note that $f^{-1}(D(c,r))$ is connected for sufficiently small $r>0
$,
and thus $U$ is completely invariant.

Since $f''(0)=0$ there is at least one parabolic basin $V$ different
from $U$ attached to the parabolic point $0$.
As $f$ has a completely invariant domain,
every component of the set of normality is simply connected.
Thus $V$ is simply connected.
Now $V$ must contain a singularity of $f^{-1}$.
Thus $V$ contains one of the critical values $d_+$ and $d_-$,
and in fact a corresponding critical point $\xi=(k\pi)^2-b$.
Since $f^n(\xi)\in\R\cap V$ and $f^n(\xi)\to 0$ as
$n\to\infty$, and since
$V$ is simply connected and symmetric with respect to the
real axis, we conclude that $(0,\xi]\subset V$.
Since $f((0,\infty))\subset (0,d_-]=f(\pi^2-b)$
we conclude that the positive real axis is contained in $V$.

We now show that $V$ is completely invariant.
Suppose that $W$ is a component of $f^{-1}(V)$ with $W\neq V$.
Since $W$ contains no critical points of $f$, and $V$ contains no
asymptotic values,
there exists a branch $\varphi$ of $f^{-1}$
which maps $V$ to $W$. This functions $\varphi$ can be continued
analytically to any point in $\bC\setminus \{ c\}$.
By the monodromy theorem, $\varphi$ extends to a
a meromorphic function from $\bC\setminus\{c\}$ to $\C$.
But this implies that $f$ is univalent, a contradiction.

It follows from part 
(iv) of our Theorem
that the Julia set of $f$ is a Jordan curve.
On the other hand we note that if $w=u+iv$ with
$|v|< T$, then
$(\Im (w^2))^2=(2uv)^2\leq 4T^2u^2< 4T^2(u^2-v^2)+4T^4
=4T^2 \Re (w^2) +4T^4$.
It follows that if
$4T^2 \Re z\leq (\Im z)^2 -4T^4$, then
$|\cos\sqrt{z}|\geq \sinh T$ and thus
$z\in U$, if
$T$ is large enough.
Thus the Julia set of $f$ is contained in the domain
$\{z\in\C: 4T^2 \Re z> (\Im z)^2 -4T^4\}$ if $T$ is large enough.
This implies that it is not a quasicircle.

\medskip

\noindent
{\bf Remark.}
It seems that $g''$ has only one negative zero.
But since we have not proved this, we have
just defined $b$ to
be the smallest zero of $g''$.
The values $a$ and $b$ are related by
$$
a= \frac{\sqrt{b}\cos\sqrt{b}-\sin\sqrt{b}}
{\sqrt{b}+\sqrt{b}\sin^2\sqrt{b}-\sin\sqrt{b}\cos\sqrt{b}}
$$
For example, if $b=-1$, then
$a=0.16763487\dots$, $g(b)=0.764166\dots$ and
$1/g'(b)=16.083479\dots$ so that
$$f(z)  =16.083479\left(\frac{1}{1+0.16763487\cos\sqrt{z-1}}-0.764166
\right).$$

\begin{figure}[htb]
\begin{center}
\includegraphics{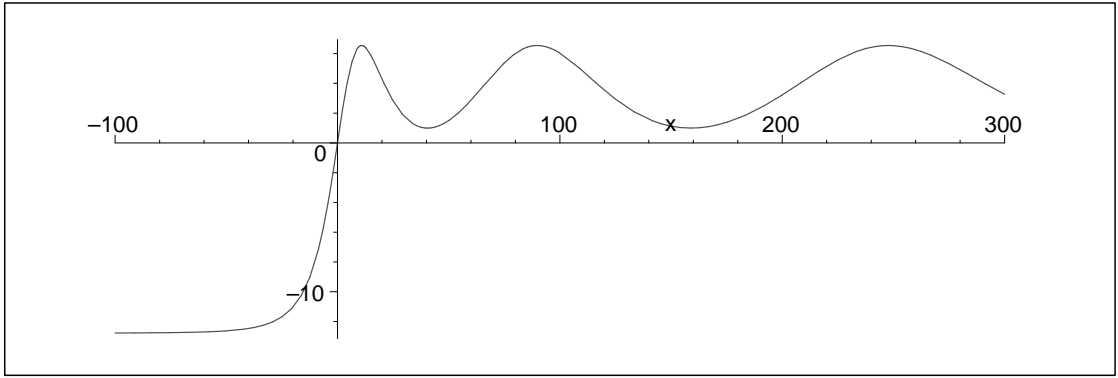}
\caption{The graph of the function $f$ from Example 1 with $b=-1$.}
\end{center}
\end{figure}

\begin{figure}[htb]
\begin{center}
\includegraphics{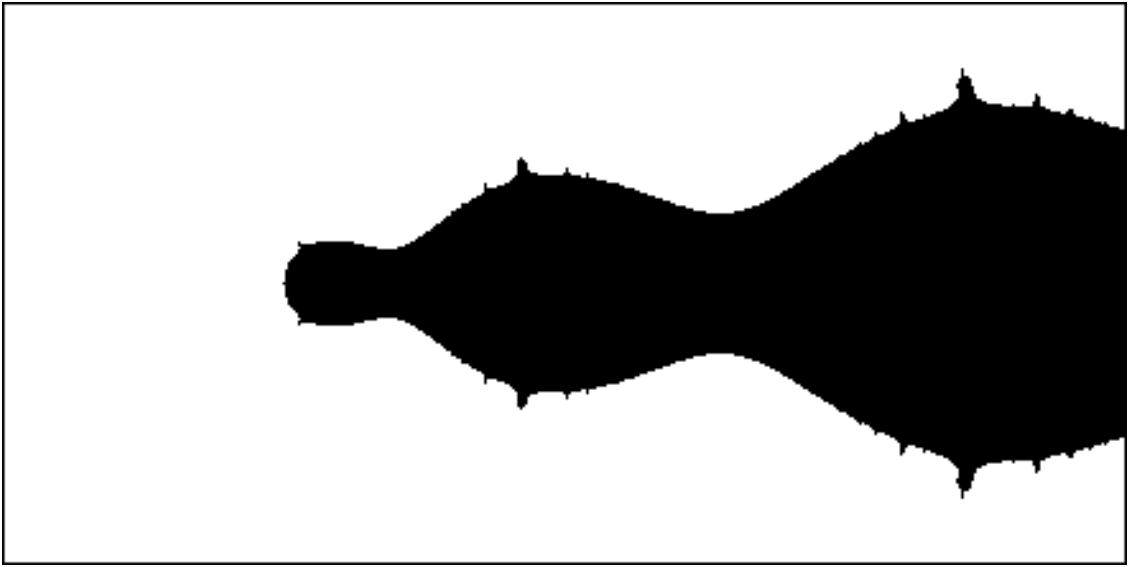}
\caption{The parabolic basin of
the function $f$ from Example 1 with $b=-1$ which contains the
positive real axis is shown in black.
The range shown is $-100<\Re z < 300, |\Im z|< 100$.
}
\end{center}
\end{figure}

\medskip

\noindent
{\bf Example 2.}
Let $g$ and $f$ be as in Example 1 and let $\alpha>1$.
Then there exists $\alpha_0=\alpha_0(a)>1$ such that
if $1<\alpha<\alpha_0$, then
$f_\alpha(z)=\alpha f(z)$ has two completely invariant
attracting basins.

\medskip
{\em Verification.}
It is not difficult to see that if
$\alpha$ is sufficiently close to $1$,
then $f_\alpha$ does
indeed have two attracting fixed points $\xi_+>0$ and
and $\xi_-<0$, with $\xi_\pm\to 0$ as $\alpha \to 1$.
The verification
that their immediate attracting basins are completely
invariant is analogous to that in Example 1.

\medskip

\noindent
{\bf Remark.}
We consider  again the case $b=-1$.
Then $f_\alpha$ has the form
$$f_\alpha(z)  =\beta \left(\frac{1}{1+0.16763487\cos\sqrt{z-1}}-0.764
166
\right)$$
with $\beta>16.083479\dots$.
For $\beta=26.712615\dots$ the positive attracting fixed point
coincides with the critical point $1+\pi^2$ and thus is 
superattracting.
\begin{figure}[htb]
\begin{center}
\includegraphics{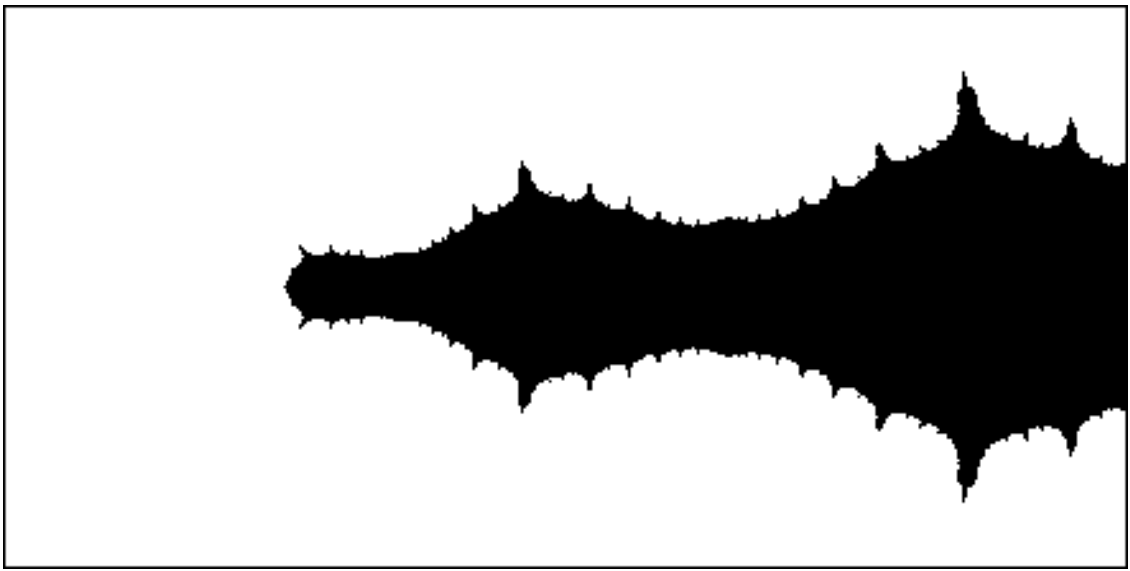}
\caption{The black region is the superattracting basin of the 
function $f$ from Example 2, with $b=-1$ and $\alpha$ chosen
such that $1+\pi^2$ is a superattracting fixed point.
The range shown is $-100<\Re z < 300, |\Im z|< 100$.
}
\end{center}
\end{figure}

\medskip

\noindent
{\bf Example 3.}
Let 
$$g(z)=\sum_{k=0}^\infty\frac{1}{a_k-z},\quad\mbox{where}\quad 0<a_0<a_1<\ldots,
\quad\;\sum_{k=0}^\infty\frac{1}{a_k}<\infty.$$
Both the upper and lower half-plane are $g$-invariant, and $g(x)\to 0$
as $x\to-\infty$ along the negative ray, so $0$ is an asymptotic value.
Evidently, the second derivative $g''$ changes sign on $(a_0,a_1)$, so there
exists $c\in(a_0,a_1)$ such that $g''(c)=0$.
Then the function
$$f(z)=\frac{g(z+c)-g(c)}{g'(c)}=z+O(z^3),\quad z\to 0$$
has a neutral fixed point with two petals at $0$.
It follows that the Julia set $J(f)$ coincides with the real line,
and thus $0\in J(f)$.

To get an example $f$ where
the upper and lower half-plane are superattracting basins,
we note that 
$g$ can be chosen such that
$g'$ has a non-real zero $\tau$, and with
$a=\Im \tau/\Im g(\tau)$ and $b=\Re \tau -a\Re g(\tau)$ the function
$f(z)=ag(z)+b$ satisfies $f(\tau)=\tau$ and $f'(\tau)=0$, as well as
$f(\overline{\tau})=\overline{\tau}$ and $f'(\overline{\tau})=0$.

\medskip

\noindent
{\bf Example 4.}
For $a=-3.7488381-1.3843391i$ the function
$f(z)=a \tan z/\tan a$ has fixed points 
$\pm a$ of multiplier $1$. The Julia set
is a Jordan curve by our theorem, but clearly
not a quasicircle.

\begin{figure}[htb]
\begin{center}
\includegraphics{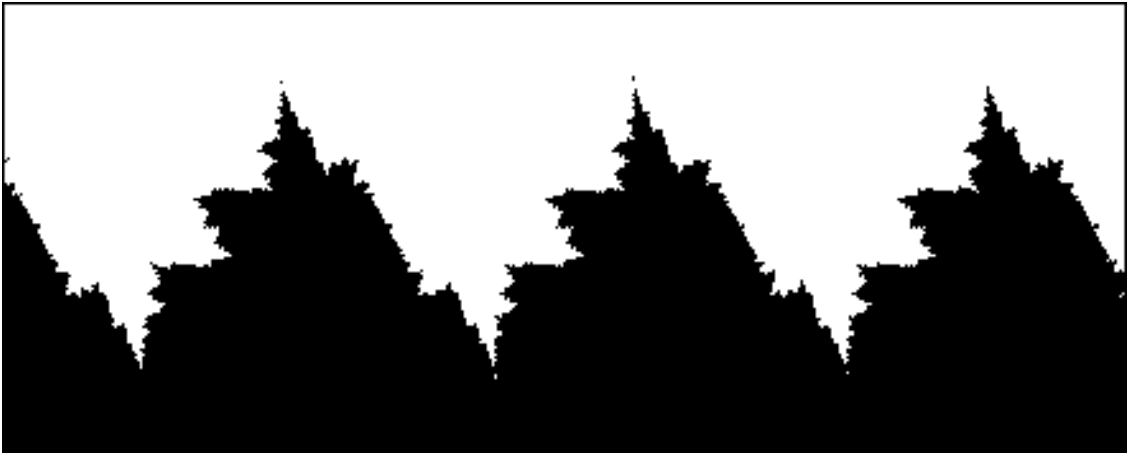}
\caption{
The parabolic basins of the function from Example 4.
The range shown is $|\Re z| < 5, |\Im z|< 2$.
}
\end{center}
\end{figure}

\medskip

\noindent
{\bf Example 5.}
For $a=1/(1-\tanh^2 1)=2.3810978$ the function
$$f(z)=a\tan z -a\tan i +i$$
has the fixed point $i$ of multiplier $1$ and the
attracting fixed point $-3.1864112i$.
Again the Julia set
is a Jordan curve, but
not a quasicircle.

\begin{figure}[htb]
\begin{center}
\includegraphics{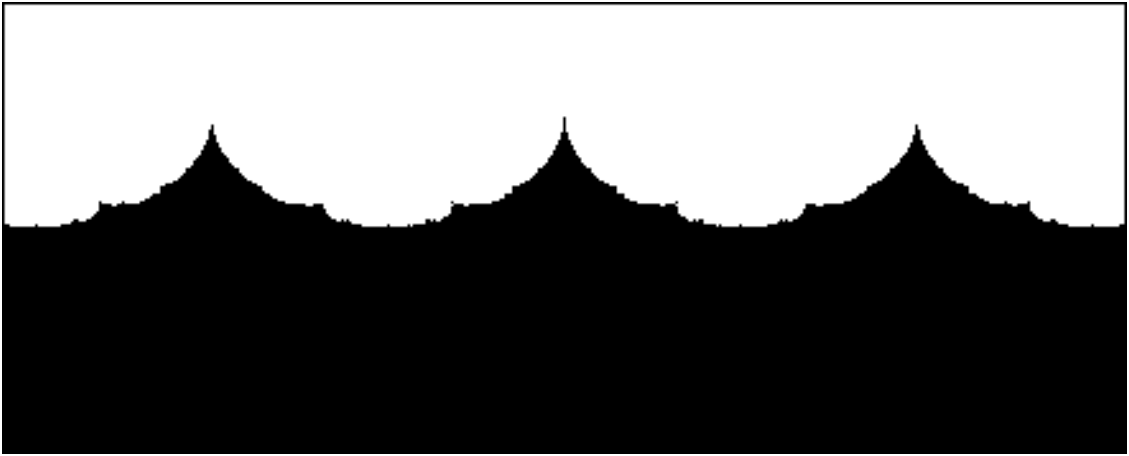}
\caption{
The parabolic basin of the function from Example 5
is shown in white, the
attracting one in black.
The range shown is $|\Re z| < 5, |\Im z|< 2$.
}
\end{center}
\end{figure}

\medskip

\noindent
{\bf Example 6.} Our final example has two completely invariant 
half-planes, but unlike $\tan z$, it has no asymptotic values.
Another feature of this example is that it has minimal possible
growth among the functions of class $S$, namely
\begin{equation}\label{gro}
T(r,f)=O((\log r)^2),\quad r\to\infty,
\end{equation}
where $T$ is the Nevanlinna characteristic. Langley
\cite{Langley,Langley2} proved that meromorphic functions
with the property $T(r,f)=o((\log r)^2)$ have infinitely
many singular values.

Let $h$ be the branch of the arccosine which
maps the 4-th quadrant $Q_4=\{ z:\Re z>0,\Im z<0\}$ onto
the half-strip $H=\{ z:\Re z\in (0,\pi/2),\Im z>0\}$.

Let $g$ be the conformal map of a
rectangle 
$R=\{ z:\Re z\in (0,\pi/2),\Im z\in (0,a)\}$
with $a>0$ onto 
$Q_4$,
such that $g(\pi/2)=0$ and $g(\pi/2+ia)=\infty$, and 
$g(ia)>g(0)>0$. By the Reflection Principle, $g$ has an analytic
continuation to the half-strip $H$ and maps this half-strip
into the left half-plane. It is easy to see that
$g$ is an elliptic function. 

The composite function $f=g\circ h$
maps the
positive ray into 
itself,
and applying the reflection
again we conclude that it maps the right half-plane into 
itself.
The boundary values on the imaginary axis
belong to the imaginary axis, so by another reflection
$f$ extends to a meromorphic function in the plane.
We see that
both the right and left half-plane are
completely invariant.

The function $f$ has 4 critical values, $\pm g(ia)$ and $\pm g(0)$,
two in the right half-plane and two in the left half-plane.

To estimate the growth of $f$ is it enough to notice that
$\arccos z=i\log z+O(1)$ as $z\to\infty$ in the lower half-plane
and in the upper half-plane. Taking into account that $g$ is
an elliptic function we obtain 
(\ref{gro}).

Our function $f$ satisfies the differential equation
$$(1-z^2)(f')^2=c(f^2-p^2)(f^2-q^2),$$
where $p=g(ia)$, $q=g(0)$ and $c$ is a real constant. 

A similar differential 
equation was considered by Bank and Kaufman~\cite{BK};
see also~\cite{Laine,Langley3}.

\begin{ack}
We thank the referee for useful comments.
\end{ack}

\end{document}